\documentclass[reqno]{amsart}
\usepackage{amsfonts}
\usepackage{bezier}
\usepackage{euscript}
\usepackage{amsfonts,amssymb,amsmath,amsbsy,amsthm,amscd}
\usepackage[english, main=english]{babel}
\usepackage{color}
\usepackage{graphics, graphicx}
\usepackage{geometry}[margin = 2in]

\newcommand{\No}{ No }

\def\dfrac#1#2{\displaystyle{#1\over #2}}

\usepackage{comment}
\usepackage{xcolor}
\usepackage[]{hyperref}
\newtheorem{theorem}{Theorem}
\begin{document}

\title[Non-Uniqueness in the Heston Model]{On Non-Uniqueness of Solutions to Degenerate Parabolic Equations in the Context of Option Pricing in the Heston Model}
\author{Ruslan~R.~Boyko} 
\address{ Mathematics and Mechanics Department, Lomonosov Moscow State University, Leninskie Gory,
Moscow, 119991,
Russian Federation}
\email{ruslan.boiko@math.msu.ru}

%\footnote{{\it Boyko Ruslan Rinatovich} --- student, Lomonosov Moscow State University, Faculty of Mechanics and Mathematics, Chair of Differential Equations, e-mail: ruslan.boiko@math.msu.ru.}}

\begin{abstract}
In \cite{HestonLoewensteinWillard2007} it was shown that the price of call options in the Heston model is determined in a non-unique way. In this paper, this problem is analyzed from the point of view of the existing mathematical theory of uniqueness classes for degenerate parabolic equations. For the special case of degeneracy, a new example is constructed demonstrating the accuracy of the uniqueness theorem for a solution in the class of functions of sublinear growth at infinity.
\end{abstract}



\maketitle
\section*{Introduction} 
In modern financial mathematics, one of the fundamental problems is finding the ``fair'' price of various derivative financial instruments, especially options, which are used to construct volatility surfaces necessary for standard calibration methods of various derivative pricing models. The uniqueness of the price guarantees the absence of arbitrage opportunities in financial markets. There are many probabilistic approaches to studying the question of price uniqueness in pricing models; however, a definitive and most complete solution to this important problem does not exist.

The problem of existence and uniqueness of option prices in the Heston model is related to the behavior of diffusion coefficients, which can degenerate when approaching the domain's boundaries. This work will explain the phenomenon of non-uniqueness in the option evaluation problem for the Heston model from the perspective of the theory of boundary value problems for degenerate second-order partial differential equations. Namely, it will be shown that non-uniqueness can arise for two reasons: the absence of necessary boundary conditions, or due to the solution leaving the Tikhonov--T\"acklind class.

\medskip
\section
{Option price in the Heston model}
The behavior of the asset in the Heston stochastic volatility model \cite{Heston1993} is described by the system of stochastic differential equations
\begin{equation}
\begin{cases}
        \label{eq:Heston1}
        d S_t = (r - q) S_t dt + \sqrt{v_t} S_t dW_t^1 \\
        d v_t =  \kappa (\theta - v_t) dt + \sigma \sqrt{v_t} dW_t^2,
\end{cases}
\end{equation}
where the volatility process $v_t$ is a Cox-Ingersoll-Ross (CIR) process \cite{CoxIngersollRoss1985}, $\theta>0$ is the long-term mean volatility, $\kappa>0$ is the rate of mean reversion, $\sigma$ is the volatility of volatility, the Brownian motions $W_t^1,$ $W_t^2$ are correlated with parameter $|\rho| \leq 1$, i.e., $dW_t^1 dW_t^2 = \rho dt$. Additionally, it is assumed that the model is written under the equivalent risk-neutral measure with drift $\mu = r - q$, where $r \geq 0$ is the interest rate and $q \geq 0$ is the dividend yield.

The option price with payoff function $\Phi(\cdot)$ in the considered model can be defined as the solution to the following boundary value problem in the domain $S, v \geq 0, \text{ and } 0 \leq t \leq T \leq \infty$ with a given terminal condition \cite{Heston1993}~\cite{GatheralTaleb2006}:
\begin{align}
\label{eq:Heston_fund}
    &\frac{v S^2}{2} V_{SS} + \rho \sigma v S V_{Sv} + \frac{v \sigma^2}{2} V_{vv} + (r - q) S V_S + (\kappa (\theta - v) - \lambda v) V_v - r V + V_t = 0,\\
\label{eq:Heston_fund_cond}
   & V(S, v, T) = \Phi(S, v).
\end{align}

\medskip
\subsection{Correct setting of boundary conditions}

Equation~\eqref{eq:Heston_fund} belongs to the class of parabolic equations degenerate on the boundary of the form
\begin{equation}\label{eq:diff_oper}
    \mathcal{L} u = \sum_{i,j=1}^{m} a_{ij}(x) u_{x_i x_j}  + \sum_{l=1}^{m} b_{l}(x) u_{x_l} + c(x) u = f(x),
\end{equation}
defined in an open domain $\mathcal{D} \subset \mathbb{R}^m$ with regular boundary $\partial \mathcal{D}$. Assume that the matrix $A(x) = \big\{a_{ij}(x)\big\}$ is symmetric, positive semidefinite (i.e., $ \langle A \xi, \xi\rangle \geq 0$ for any unit vector $ \xi$), and has degeneration points on the boundary $\partial \mathcal{D}$, with the coefficient $c(x) < 0$. A general theory of weak solutions for such equations was constructed in \cite{Fichera1963}~\cite{OleinikRadkevich1971}. To investigate the correctness of the boundary condition formulation, we introduce the \emph{Fichera function}
\[
    H(x) = \sum_{i=1}^{m} \Big[b_i(x) - \sum_{j=1}^{m} \big(a_{i j}(x)\big)_{x_j}\Big] n_i,
\]
constructed for points of the set $\Sigma^0 = \big\{x\in \partial D \mid \langle A(x) \nu, \nu \rangle = 0 \big\} \subset \partial D$, where $\nu$ denotes the unit normal to the boundary $\partial \mathcal{D}$ at point $x$.

The Fichera function for problem \eqref{eq:Heston_fund},~\eqref{eq:Heston_fund_cond} has the form
\[
    H(S, v, t) = \left((r-q) S - \Big(v S + \frac{\rho \sigma S}{2}\Big)\right) n_1 + \left(\kappa (\theta - v) - \lambda v - \Big(\frac{\rho \sigma v}{2} + \frac{\sigma^2}{2}\Big)\right) n_2 + n_3,
\]
where $\mathbf{n} = (n_1, n_2, n_3)$ defines the inward normal in coordinates $(S, v, t)$ to the boundary of the domain.

Let us introduce the domains associated with the sign of the Fichera function
\begin{align*}
    \Sigma_0 &= \big\{x \in \Sigma^0 \mid H(x) = 0\big\}, &&\Sigma_2 = \big\{x \in \Sigma^0 \mid H(x) < 0\big\}, \\
    \Sigma_1 &= \big\{x \in \Sigma^0 \mid H(x) > 0\big\}, &&\Sigma_3 = \big\{\partial D - \Sigma^0\big\}.
\end{align*}
For \eqref{eq:diff_oper}, we pose the boundary value problem $u = g \text{ on } \Sigma_2 \cup \Sigma_3$. In \cite{Fichera1963}, a theorem on the existence of a weak solution applicable to our problem was proved. Namely, assume that for the boundary value problem
\begin{quote}
    \it the following conditions hold: $c(y) \leq c_0 < 0$ in $D$, $f$ is a bounded measurable function on $D$, $g$ is a bounded measurable function on $\Sigma_{2} \cup \Sigma_{3}$. Then there exists a weak solution to the posed boundary value problem.
\end{quote}
Thus, for the correct formulation of the problem, it is necessary to impose a boundary condition on the domain $\Sigma_2 \cup \Sigma_3$, i.e., at points of non-degeneracy and at points where the Fichera function is negative.

According to Fichera's theory \cite{Fichera1963}~\cite{Meyer2015}, for the well-posedness of the boundary value problem it is necessary to impose a boundary condition at the terminal moment $t = T$. Furthermore, if the model parameters do not satisfy the inequality
\begin{equation}\label{Feller}
  k \theta \geqslant \frac{\sigma^2}{2},
\end{equation}
then a boundary condition must also be imposed on the boundary $v = 0$ (the moment the volatility trajectory reaches the zero level). This condition is commonly referred to in probability theory as the \emph{Feller condition}.

Thus, under certain conditions on the model parameters, the non-uniqueness of the solution in problem \eqref{eq:Heston_fund}, \eqref{eq:Heston_fund_cond} occurs due to a missing boundary condition.

Studying the sign of the Fichera function shows that on the boundaries $S=0$, $S = +\infty$, and $v = +\infty$, imposing a boundary condition is not required. Note that the analysis of the necessity of imposing boundary conditions for the option price equation in the Heston model was carried out in \cite{Meyer2015}.

\medskip

\section{ The uniqueness problem in option pricing in the Heston model under the Feller condition}

In \cite{HestonLoewensteinWillard2007}, an Example 3 showing that the solution to the option pricing problem may be non-unique. For this, a special case of model \eqref{eq:Heston1} was considered:
\begin{equation*}
\begin{cases}
        \label{eq:Heston2}
        d S_t = r S_t dt + \sqrt{\nu_t} S_t dZ^{\mathbb{Q}}_t\\
        d \nu_t = \sigma^2 dt + \sigma \sqrt{\nu_t} dZ^{\mathbb{Q}}_t,
\end{cases}
\end{equation*}
where the parameters $\sigma > 0, r \geq 0$ correspond to the characteristics of the volatility process and the drift of the price process, respectively, and $Z^{\mathbb{Q}}_t$ is a Wiener process under the risk-neutral measure~$\mathbb{Q}$.

Equation \eqref{eq:Heston_fund} becomes
\begin{align}\label{eq:H1}
    \frac{v S^2}{2} V_{SS} + \sigma v S V_{Sv} + \frac{v \sigma^2}{2} V_{vv} + r S V_S + \sigma^2 V_v - r V + V_t &= 0.
\end{align}

Note that from the standpoint  of Fichera function theory, this equation, as in the general Heston model, requires a condition to be imposed only at the terminal time $t = T$ (the investigation reduces to examining the sign of the Fichera function on the boundaries). Here, the "Feller condition," which previously arose on the boundary $v = 0$, is replaced in the considered problem by the condition $\sigma^2 - \frac{\sigma^2}{2} = \frac{\sigma^2}{2} \geqslant 0$, i.e., imposing a condition on the boundary $v = 0$ is not required.

In \cite{HestonLoewensteinWillard2007}, the differential equation
\begin{equation}
\label{eq:Pi_add}
    \frac{v \sigma^2}{2} \Pi_{vv} + \sigma^2 \Pi_v + \Pi_t - r \Pi = 0
\end{equation}
with final condition $\Pi(v, T) = 0$ was considered to prove the non-uniqueness of the solution to problem \eqref{eq:H1}, \eqref{eq:Heston_fund_cond}. This work shows that the function
\begin{equation}
\label{eq:Heston_solution}
    \Pi(v, t) = \frac{1}{v} e^{-r (T-t) - \frac{2 v}{\sigma^2 (T-t)}}
\end{equation}
provides a nontrivial solution to this problem. Consequently, if $V_1 (S, v, t)$ denotes the ``standard'' solution of problem \eqref{eq:H1}, \eqref{eq:Heston_fund_cond} (e.g., \cite{Heston1993}~\cite{GatheralTaleb2006}), then $V_2(S, v, t) = V_1(S, v, t) + \Pi(v, t)$ also solves the original problem.

Note that solution \eqref{eq:Heston_solution} is unbounded as $v\to 0$, and it is natural to assume that to single out a unique solution, restrictions on the behavior of the solution as $v\to 0$ must be imposed. Below it will be shown how the reason for this kind of non-uniqueness can be explained in terms of classical results of the theory of parabolic equations.

Furthermore, questions concerning the growth restriction of the solution also arise as $S\to 0$, $v\to +\infty$, $S\to +\infty$. We will also obtain growth conditions for the solution that ensure uniqueness when approaching finite or infinite boundaries.

\medskip
\section{Classes of uniqueness in the problem of option valuation in the Heston model}\label{S3}
\subsection{ Tikhonov--T\"acklind classes  ($S\to\infty$, $v\to\infty$).}

The question of identifying uniqueness classes for the Cauchy problem for degenerate parabolic equations was first addressed by E.~Holmgren in his 1924 work. The exact uniqueness class for the solution of the Cauchy problem in the context of the heat equation was obtained by A.N.~Tikhonov in 1935 in \cite{Tikhonov1935}. It was established that the solution $u(x,t)$ of the heat equation $u_t = \Delta u$ with zero initial condition $u(x, 0) = 0$ is unique in the class of functions
\[
    |u(x,t)| \leq B \exp\{\beta|x|^2\}, \quad \beta > 0,\, x \in \mathbb{R}^n,\, t \in [0,T],
\]
and Tikhonov constructed an example of a nonzero solution to this Cauchy problem that belonged to a wider class. In 1936, S.~T\"acklind in \cite{Tacklind1936} refined Tikhonov's results by showing that the solution of the Cauchy problem is unique in the class
\[
    |u(x,t)| \leq B \exp\{|x| h(x)\}, \quad x \in {\mathbb R}, \, t \in [0,T],
\]
where $h(\cdot)$ is a non-decreasing nonnegative function with a divergent integral $\int_1^{+\infty} \frac{ds}{h(x)}$ (this function is commonly called the \emph{T\"acklind function}). It was also established that if the latter condition is violated, the solution ceases to be unique.

The works of Tikhonov and T\"acklind were successively generalized to broader classes of equations. In particular, Kamynin and Khimchenko in \cite{KamyninKhimchenko1979}~\cite{KamyninKhimchenko1981}~\cite{Kamynin1984} extended the theory to general second-order parabolic equations with nonnegative characteristic form and unbounded coefficients. They proposed a generalized Tikhonov--T\"acklind class
\[
    |u(x,t)| \leq C \exp\big\{ G(|x|) h(G(|x|)) \big\}, \quad  G(s) = \int_0^s \frac{dx}{g(x)},
\]
where $g\colon [0, +\infty) \to [1, +\infty)$ is a non-decreasing function of class $C^1[0, +\infty)$, $g(s) \equiv 1$ for $s \in [0, 1]$ and $\int_0^{+\infty} \frac{ds}{g(s)} = \infty$.

In \cite{KamyninKhimchenko1981}, Theorem 1 was proved, which we formulate in a special case that allows its application to our situation. Namely, assume that
\begin{quote}\it
     the coefficients of the differential operator
     \[
         \mathcal{L} = \sum_{i,j=1}^{m} a_{ij}(x, t) \partial_{x_i x_j}  + \sum_{l=1}^{m} b_{l}(x, t) \partial_{x_l} + c(x, t) - \partial_t ,
     \]
     in the strip $\Pi(T) = \{x \in \mathbb{R},\, 0 \leqslant t \leqslant T\}$ satisfy the conditions
    \begin{align*}
       & 0 \leq A(x,t; \xi) = \sum_{i, j = 1}^{n} a_{ij}(x, t) \xi_i \xi_j \leq k(|x|), \quad \forall |\xi| = 1, \\
        &0 \leq a(x,t) = \left( \sum_{i=1}^{n} b_i^2(x, t)\right)^{1/2} \leq g(|x|) \varphi(G(|x|)), \quad
        c(x,t) \leq 0,
    \end{align*}
    where $k(s) = \min\{g^2(s), s g(s)\}$ for $s \geq 1$, and the function $u(x,t) \in C(\Pi(T)) \cap C_{x,t}^{2,1}(\Pi(T))$ belongs to the generalized Tikhonov--T\"acklind class. Then, if   $\mathcal{L}u = 0$ and $u(x,0) = 0$, then $u(x,t) \equiv 0$ in the entire strip $\Pi(T)$.
\end{quote}

To find the Tikhonov--T\"acklind uniqueness class for equation \eqref{eq:H1} as $S \to +\infty$ and \mbox{$v \to +\infty$}, we perform the substitution $x = \ln S$. Then equation \eqref{eq:H1} is rewritten as
\begin{equation}\label{HestonX}
    \frac{v}{2}V_{xx} + \sigma v V_{xv} + \frac{v\sigma^2}{2}V_{vv} + \left(r - \frac{v}{2}\right)V_x + \sigma^2 V_v - rV + V_t = 0.
\end{equation}
Let us estimate the quadratic and linear coefficients:
\begin{align*}
    0\le A(x,v,t;\xi) &= \frac{v}{2}\xi_1^2 + \sigma v\xi_1\xi_2 + \frac{v\sigma^2}{2}\xi_2^2 \leq C_1 v, \\
    0\le a(x, v, t) &= \left( \left(r - \frac{v}{2}\right)^2 + (\sigma^2)^2 \right)^{1/2} \leq C_2 v,
\end{align*}
where $C_1$, $C_2$ are constants.
Choose
\[
    g(s) = \begin{cases}
        1, &s \in [0, 1 - \varepsilon] \\
        \text{smooth}, &s \in [1-\varepsilon, 1+\varepsilon)\\
        s, &s \geqslant 1 + \varepsilon
    \end{cases} \qquad G(s) = \int_0^s \frac{dx}{g(x)} \sim C \ln(s), \, s\to +\infty,
\]
which satisfy the necessary conditions given in \cite{KamyninKhimchenko1979}~\cite{KamyninKhimchenko1981}.
Thus, from Theorem 1 in \cite{KamyninKhimchenko1979} follows that the uniqueness class as $S \to  \infty$  and $v \to \infty$ has the form:
\begin{align}\label{V}
    |V(S, v, t)| \leqslant C \exp \left[ \ln\left(\sqrt{\ln^2(S) + v^2}\right) \cdot h \left( \ln\left(\sqrt{\ln^2(S) + v^2}\right) \right) \right].
\end{align}
\subsection{ Uniqueness class for  \eqref{eq:H1} as $S\to 0$.}

Note that the substitution $x=\ln S$ also allows us to investigate the behavior of the solution as $S\to 0$, since it reduces the problem to studying the uniqueness class for equation \eqref{HestonX} as $x\to -\infty$. This equation does not degenerate in the variable $x$, so the uniqueness class coincides with the well-known uniqueness class for the heat equation.
Performing the inverse substitution, we see that condition \eqref{V} also singles out the uniqueness class for equation \eqref{eq:H1} as $S\to 0$.

\medskip
\subsection{ Uniqueness class for  \eqref{eq:H1} as $v\to 0$.}

The question of uniqueness classes in this case is resolved somewhat differently. We could, of course, make the substitution $w=\frac{1}{v}$ and study the behavior of the solutions of the resulting equation as $w\to\infty$, but the theory used above turns out to be inapplicable in this case. Therefore, we will proceed differently and study the behavior of the solution of equation
\eqref{eq:Pi_add} as $v\to 0$ within the framework of the theory of degenerate partial differential equations, similar to what was done in~\cite{CoxHobson2005}.

Note that equation~\eqref{eq:Pi_add} reduces to the Feller equation~\cite{Feller1951}
\begin{equation}
\label{eq:Feller}
    u_t = (a x u)_{xx} - ((b x + c) u)_{x},
\end{equation}
where $a, b, c$ are some constants. Namely, performing the substitution $x = \frac{4}{\sigma^2} v$ and $\tau = T - t$ (assuming $r = 0$) we obtain
\[
    \Pi_{\tau} = 2 x \Pi_{xx} + 4 \Pi_x,
\]
i.e., $a = 2, b = 0$ and $c = 4 = \frac{3 - 2 \alpha}{1 - \alpha}, \alpha = \frac{1}{2}$.

In turn, the substitution $y = \frac{1}{x}$ leads us to the equation
\begin{equation}
\label{eq:Pi1}
    \Pi_{\tau} = 2 y^3 \Pi_{yy},
\end{equation}
with parameter $\alpha = \frac{3}{2}$. Equation \eqref{eq:Pi1} belongs to the class of equations with power-law degeneracy of the form
\begin{equation}
\label{eq:Pi2}
    \Pi_{\tau}= a y^{2\alpha}
    \Pi_{yy},\quad \alpha>0,\, a>0,
\end{equation}
which are closely related to stochastic processes with constant elasticity (CEV) \cite{Hull2019}. It is known that the change of variables
\[
    x = \frac{1}{2 a (1-\alpha)^2} y^{2(1-\alpha)}
\]
reduces them to the Feller equation \eqref{eq:Feller} \cite{Hull2019}. For equation \eqref{eq:Pi2} in the case $\alpha>1$, there also exists a non-uniqueness problem, which was investigated in \cite{LadykovaRozanova2025} from the perspective of the theory of degenerate partial differential equations. Thus, \eqref{eq:Pi1} is a special case of such a problem.

Let us list the results concerning the growth rate of the solution as $y\to\infty$ that ensure its uniqueness:

1. In \cite{EkstromTysk2009} (Theorem 4.8) it is shown that for $\alpha \geq 0$,
%Theorem 4.3 [13] gives a sufficient condition for the uniqueness of the solution in the posed problem:
\begin{quote}\it
    if the initial condition $\Pi(x, 0) = g(x)$ in \eqref{eq:Pi_add} has (strictly) sublinear growth at infinity, then the solution to problem \eqref{eq:Pi_add} with this terminal condition is unique.
\end{quote}

2. In \cite{LadykovaRozanova2025} it is shown that for $0\leq \alpha \leq 1$ the standard approach related to the construction of exact Tikhonov--T\"acklind uniqueness classes \cite{KamyninKhimchenko1981}~\cite{Kamynin1984} is applicable; functions from this class can grow faster than a linear function, i.e., the necessary condition given in \cite{EkstromTysk2009} is not sufficient.

3. For $\alpha>1$, the theory of \cite{KamyninKhimchenko1979}~\cite{KamyninKhimchenko1981} is not applicable and, apparently, it is unknown whether the sublinear growth condition is sufficient for identifying uniqueness classes. At the same time, an example of a nontrivial solution in the case of zero initial data for $\alpha = 2$ is known (see \cite{CoxHobson2005})
\begin{equation}\label{Example}
   {\Pi}_2 (y, t) = y \left(1 - 2 \Phi \left(- \frac{1}{\sigma y \sqrt{t}}\right)\right),
\end{equation}
   where $\Phi(z) = \dfrac{1}{\sqrt{2 \pi}} \int_{-\infty}^{z} e^{-u^2/2} \, du$ is the Laplace function. This solution has linear growth as $y\to+\infty$ and serves as an example of non-uniqueness of the solution for the given parameter value. Thus, for $\alpha = 2$, Theorem 4.3 of \cite{EkstromTysk2009} gives both a necessary and sufficient condition for the uniqueness of the solution.

Let us construct a similar example of non-uniqueness for $\alpha = 3/2$ based on solution \eqref{eq:Heston_solution}. For this, it suffices to note that \eqref{eq:Heston_solution} under the substitution will become
\begin{equation}
\label{eq:Heston_solution_inf}
    \Pi(y, \tau) = \frac{4 y}{\sigma^2} e^{- \frac{1}{2 y \tau}}, \quad y>0, \quad \tau>0.
\end{equation}
This function is a nontrivial solution of \eqref{eq:Pi2}, vanishes at $\tau=0$, and grows linearly as $y\to+\infty$.

Thus, the uniqueness class for solutions in problem \eqref{eq:H1}, \eqref{eq:Heston_fund_cond} as $v \to +0$ consists of functions having sublinear growth at infinity (after the coordinate change $y = \frac{4}{\sigma^2 v}$), i.e., the solution has a singularity weaker than $1/v$, in particular, integrable. Outside this class, the solution of the considered equation is, generally speaking, non-unique.

\medskip
\subsection{Main theorem}

Let us summarize the results obtained in Section \ref{S3} on uniqueness classes for problem \eqref{eq:H1}, \eqref{eq:Heston_fund_cond}:

\smallskip
\begin{theorem}
  Assume that the parameters of model \eqref{eq:Heston1} are such that the Feller condition \eqref{Feller} is satisfied. Then the classical solution to problem \eqref{eq:H1}, \eqref{eq:Heston_fund_cond} is unique in the class of functions having, as $S\to +\infty$, $v\to +\infty$, $S\to 0$, growth determined by condition \eqref{V}, and belonging to the class $o\left(\frac{1}{v}\right)$ as  $v\to 0$.
\end{theorem}

\medskip
\section{ Conclusion and discussion}

 This work contains the following main results.

 1. The uniqueness class in the option pricing problem in the Heston model has been found. In particular, it is shown that it consists of functions having a singularity weaker than $1/v$ as $v\to+0$. This shows the reason for non-uniqueness in the example constructed in \cite{HestonLoewensteinWillard2007}.

 2. Example \eqref{eq:Heston_solution_inf} is constructed, showing that the requirement of sublinear growth at infinity is not only necessary but also a sufficient condition for the uniqueness of the solution to the parabolic equation with power-law degeneracy \eqref{eq:Pi2} for $\alpha=3/2$.
  In fact, in addition to the known example \eqref{Example} for $\alpha=2$, we have constructed a second example of this kind.

A natural question arises: is it possible, by reduction to the Feller equation, to construct a similar non-uniqueness example for the remaining $\alpha>1$ and thereby show that the sublinear growth condition at infinity is precise for identifying the uniqueness class for solutions of equations of the form \eqref{eq:Pi2}? Such a result seems quite natural, but the corresponding examples are currently unknown. A formal obstacle to transferring the result from $\alpha=3/2$ to other $\alpha>1$ is that when reducing \eqref{eq:Pi1} to the Feller equation for $\alpha\ne 3/2$, additional lower-order terms appear.

\end{document}